\documentclass{article}
\usepackage{graphicx}
\usepackage{amsmath}
\usepackage{hyperref}
\usepackage{amsfonts}
\usepackage{amssymb}
\newtheorem{theorem}{Theorem}[section]

\newtheorem{corollary}[theorem]{Corollary}

\newtheorem{definition}[theorem]{Definition}
\newtheorem{example}[theorem]{Example}

\newtheorem{notation}[theorem]{Notation}

\newtheorem{proposition}[theorem]{Proposition}
\newtheorem{remark}[theorem]{Remark}

\numberwithin{equation}{section}
\newenvironment{proof}[1][Proof]{\noindent\textbf{#1\ } }{\rule{0.5em}{0.5em}\medskip}

\begin{document}

\title{Relative and Discrete Utility Maximising Entropy}
\author{Grzegorz Hara\'{n}czyk\\{\small Institute of Mathematics, Jagiellonian University}\\{\small Reymonta 4, 30-059 Krak\'{o}w, Poland}\\{\small grzegorz.haranczyk@im.uj.edu.pl}
\and Wojciech S\l omczy\'{n}ski\\{\small Institute of Mathematics, Jagiellonian University}\\{\small Reymonta 4, 30-059 Krak\'{o}w, Poland}\\{\small wojciech.slomczynski@im.uj.edu.pl}
\and Tomasz Zastawniak\\{\small Department of Mathematics, University of York}\\{\small Heslington, York YO10 5DD, UK}\\{\small tz506@york.ac.uk}}
\date{}
\maketitle
\begin{abstract}
\noindent The notion of utility maximising entropy ($u$-entropy) of a
probability density, which was introduced and studied in~\cite{SloZas04}, is
extended in two directions. First, the relative $u$-entropy of two probability
measures in arbitrary probability spaces is defined. Then, specialising to
discrete probability spaces, we also introduce the absolute $u$-entropy of a
probability measure. Both notions are based on the idea, borrowed from
mathematical finance, of maximising the expected utility of the terminal
wealth of an investor. Moreover, $u$-entropy is also relevant in
thermodynamics, as it can replace the standard Boltzmann-Shannon entropy in
the Second Law. If the utility function is logarithmic or isoelastic (a power
function), then the well-known notions of the Boltzmann-Shannon and R\'{e}nyi
relative entropy are recovered. We establish the principal properties of
relative and discrete $u$-entropy and discuss the links with several related
approaches in the literature.
\end{abstract}

\section{Introduction}

The notion of utility maximising entropy (or $u$-entropy for brevity) of a
density~$f$ of a probability distribution with respect to a given probability
measure~$\mu$ was introduced and studied by two of the present authors
in~\cite{SloZas04}.

The work in \cite{SloZas04} was motivated, on the one hand, by problems in
mathematical finance concerned with a trader with a concave utility
function~$u$ who wants to maximise the expected utility $\mathbb{E}_{\nu
}(u(w))$ under the true market probability measure~$\nu$ over all contingent
claims~$w$ (where~$w$ is a non-negative random variable representing the final
value of a contingent claim) whose initial value $\mathbb{E}_{\mu}(w)$ under a
pricing measure~$\mu$ is equal to the initial wealth of the trader, taken to
be~$1$ for simplicity, so that $\mathbb{E}_{\mu}(w)=1$. A value $c\in
\mathbb{R}$, called the certainty equivalent, can be assigned to each
contingent claim~$w$ so that $u(c)=\mathbb{E}_{\nu}(u(w))$. If $\nu$ is
absolutely continuous with respect to~$\mu$ with density~$f$, then the
$u$-entropy $H_{u}(f)$ is defined as the highest possible value of the
logarithm of the certainty equivalent~$c$ over all contingent claims $w$ with
initial value $\mathbb{E}_{\mu}(w)=1$.

Expected utility maximisation problems in mathematical finance have been
studied extensively, for example, in \cite{PikKar96}, \cite{AmeImkSch98},
\cite{RouKar00}, \cite{GolRus01}, \cite{BelFri02}, \cite{FriBia05},
\cite{Gun05}, \cite{GunFol07}. Some of the most general and elegant results,
which have provided much inspiration for our work, belong to Kramkov and
Schachermayer \cite{KraSch99}, \cite{KraSch03}, \cite{Sch04},
\cite{HugKraSch05}.

On the other hand, further motivation for $u$-entropy comes from
thermodynamics and statistical mechanics. When $u(x)=\ln x$, then the
$u$-entropy $H_{u}(f)$ is equal to the classical Boltzmann-Gibbs entropy
$H(f)=\mathbb{E}_{\mu}(f\ln f)$ (note the sign convention typical of
mathematical literature; the opposite sign for entropy would normally be used
in physics). The properties of the Boltzmann-Gibbs entropy $H(f)\ $and, in
particular, its role in the Second Law of thermodynamics provided a fertile
ground for generalisation to the case of $u$-entropy. As is very well known, a
physical system in state~$f$ evolves towards equilibrium whenever $H(f)$ tends
to zero. An extension of this and other properties of entropy was achieved
in~\cite{SloZas04} by replacing the Boltzmann-Gibbs entropy $H(f)$ with the
$u$-entropy $H_{u}(f)$ for an arbitrary $u$ from a broad class of utility functions.

In the present paper the concept of $u$-entropy is extended further to include
relative entropy of two arbitrary probability measures $\nu$ and~$\mu$. It
will be called the relative\ $u$-entropy and denoted by $H_{u}\left(
\nu\parallel\mu\right)  $. We also introduce the $u$-entropy $h_{u}(p)$ of a
probability measure $p$ (rather than the relative entropy of one measure with
respect to another or that of a density with respect to a given probability
measure), but to do so need to specialise to the case of a discrete
probability space, where $p$ is a probability vector.

We establish some properties of relative $u$-entropy and discrete $u$-entropy,
and study their relationships with other similar approaches in the literature.
In particular, we discuss a link with the recent work by Friedman, Huang and
Sandow \cite{FriHuaSan07}, and with a much older approach by Arimoto
\cite{Ari71}, which does not refer to utility maximisation explicitly but is
based on a similar concept. These two approaches work in the discrete case
only. Moreover, returning once again to general probability spaces, we also
establish a connection of relative $u$-entropy with Frittelli's generalised
distance between two probability measures, introduced in~\cite{Fri00a} to
solve the dual convex problem in a utility maximisation framework for asset
pricing in an incomplete market.

It will prove convenient to adopt the convention $\infty\cdot0=-\infty
\cdot0=0$ throughout this paper.

\section{Utility maximising relative entropy}

\subsection{Utility functions}

\begin{definition}
\upshape Let $u:(0,\infty)\mathbb{\rightarrow}\mathbb{R}$. We call $u$ a
\textbf{utility function} whenever $u$ satisfies the \textbf{Inada
conditions}, that is, $u$~is a strictly concave strictly increasing
continuously\ differentiable function such that
\[
u^{\prime}\left(  0\right)  :=\lim_{x\searrow0}u^{\prime}\left(  x\right)
=\infty,\quad u^{\prime}\left(  \infty\right)  :=\lim_{x\nearrow\infty
}u^{\prime}\left(  x\right)  =0\text{ .}%
\]
We shall also use the notation
\[
u\left(  0\right)  :=\lim_{x\searrow0}u\left(  x\right)  ,\quad u\left(
\infty\right)  :=\lim_{x\nearrow\infty}u\left(  x\right)  \text{ .}%
\]
\end{definition}

\begin{proposition}
The function $I:=\left(  u^{\prime}\right)  ^{-1}:(0,\infty)\rightarrow
(0,\infty)$ is strictly decreasing and satisfies
\[
I\left(  0\right)  :=\lim_{x\searrow0}I\left(  x\right)  =\infty,\quad
I\left(  \infty\right)  :=\lim_{x\nearrow\infty}I\left(  x\right)  =0\text{ .}%
\]
\end{proposition}

\begin{definition}
\upshape
Let $u:(0,\infty)\rightarrow\mathbb{R}$ be a utility function. The
\textbf{convex dual} $u^{\ast}:(0,\infty)\rightarrow\mathbb{R}$ is defined by
\begin{equation}
u^{\ast}\left(  y\right)  =\sup_{x>0}\left(  u\left(  x\right)  -yx\right)
\label{defconuti}%
\end{equation}
for any $y\in(0,\infty)$. We also put
\begin{equation}
u^{\ast}(0):=\lim_{x\searrow0}u^{\ast}(x)=u(\infty),\quad u^{\ast}\left(
\infty\right)  :=\lim_{x\nearrow\infty}u^{\ast}\left(  x\right)  =u(0).
\label{defconinf}%
\end{equation}
\end{definition}

If $\Lambda>0$ and $s=0$, we put $I\left(  \Lambda/s\right)  :=0$ and
$u^{\ast}\left(  \Lambda/s\right)  :=u\left(  0\right)  $, consistently with
the adopted notation $I\left(  \infty\right)  =0$ and $u^{\ast}\left(
\infty\right)  =u\left(  0\right)  $.

\begin{example}
\upshape
Let $\gamma\in\left(  -\infty,1\right)  $. Define $u:(0,\infty)\rightarrow
\mathbb{R}$ by
\[
u\left(  t\right)  =\left\{
\begin{array}
[c]{cl}%
\frac{1}{\gamma}\left(  t^{\gamma}-1\right)  & \text{for }t\in(0,\infty)\text{
and }\gamma\in\left(  -\infty,0\right)  \cup\left(  0,1\right) \\
\ln t & \text{for }t\in(0,\infty)\text{ and }\gamma=0
\end{array}
\right.  \text{.}%
\]
We call $u$ the \textbf{isoelastic utility of order }$\gamma$ if $\gamma\neq
0$, and the \textbf{logarithmic utility} if $\gamma=0$.
\end{example}

The following definition is due to Kramkov and Schachermayer~\cite{KraSch99}.

\begin{definition}
\upshape
The \textbf{asymptotic elasticity} of a utility function $u:(0,\infty
)\rightarrow\mathbb{R}$ is defined by
\[
\mathrm{AE}(u)=\limsup_{x\nearrow\infty}\frac{xu^{\prime}\left(  x\right)
}{u\left(  x\right)  }\text{ .}%
\]
A utility function $u$ is said to have \textbf{reasonable asymptotic
elasticity} if $\mathrm{AE}(u)<1$.
\end{definition}

Under the assumption of reasonable asymptotic elasticity, duality theory for
utility maximisation works in a similar manner as in the finite-dimensional
case. See \cite{Sch04} for equivalent formulations of this assumption and a
discussion of its economic meaning.

\subsection{Relative $u$-entropy and $u$-entropy}

\subsubsection{Definition}

\begin{notation}
\upshape
Let $\left(  \Omega,\Sigma\right)  $ be a measurable space. We denote by
$M_{1}\!\left(  \Omega,\Sigma\right)  $ the space of all probability measures
on $\left(  \Omega,\Sigma\right)  $. For any $\mu\in M_{1}\!\left(
\Omega,\Sigma\right)  $ we denote by $D\left(  \mu\right)  $ the set of all
\textbf{densities} on the probability space $\left(  \Omega,\Sigma,\mu\right)
$, that is,
\[
D\left(  \mu\right)  :=\left\{  w\in L^{1}\left(  \mu\right)  :w\geq0\text{
and }\int_{\Omega}wd\mu=1\right\}  \text{ .}%
\]
By $B\left(  \Omega,\Sigma\right)  $ we denote the set of all bounded
measurable real-valued functions on $\left(  \Omega,\Sigma\right)  $. In the
sequel we shall write simply $M_{1}$ and $B$ whenever the measurable space
$\left(  \Omega,\Sigma\right)  $ is unambiguous. For any $\mu\in M_{1}$ and
$f\in D(\mu)$ we shall write~$f\mu$ to denote the measure in~$M_{1}$ with
density~$f$ with respect to~$\mu$.
\end{notation}

\begin{definition}
\label{defreluentr}\upshape
Let $u:(0,\infty)\mathbb{\rightarrow}\mathbb{R}$ be a utility function. Let
$\left(  \Omega,\Sigma\right)  $ be a measurable space and let $\nu,\mu\in
M_{1}\left(  \Omega,\Sigma\right)  $. We put
\[
N_{u}\left(  \nu\parallel\mu\right)  :=\sup_{w\in\mathcal{A}\left(  \nu
,\mu\right)  }\int_{\Omega}u\left(  w\right)  d\nu\text{ ,}%
\]
where
\[
\mathcal{A}\left(  \nu,\mu\right)  :=\left\{  w\in D\left(  \mu\right)
:u\left(  w\right)  ^{-}\in L^{1}\left(  \nu\right)  \right\}  \text{ .}%
\]
Here $x^{-}=\max(-x,0)$ denotes the negative part of $x\in\mathbb{R}$. Note
that $\int_{\Omega}u(w)d\nu\in(-\infty,\infty]$ for each $w\in\mathcal{A}%
\left(  \nu,\mu\right)  $. We define
\[
H_{u}\left(  \nu\parallel\mu\right)  :=\ln u^{-1}\left(  N_{u}\left(
\nu\parallel\mu\right)  \right)
\]
and call it the \textbf{relative }$u$-\textbf{entropy} (or \textbf{relative
utility maximising entropy}) of~$\nu$\ with respect to~$\mu$.
\end{definition}

\begin{definition}
\label{defutmaxent}\upshape
(from \cite{SloZas04}) Let $u:(0,\infty)\rightarrow\mathbb{R}$ be a utility
function and let $\mu\in M_{1}$. For any $f\in D(\mu)$ we put
\[
N_{u}\left(  f\right)  :=\sup_{w\in\mathcal{A}\left(  f\right)  }\int_{\Omega
}u\left(  w\right)  f\,d\mu\text{ ,}%
\]
where
\[
\mathcal{A}\left(  f\right)  :=\left\{  w\in D\left(  \mu\right)  :u\left(
w\right)  ^{-}\in L^{1}\left(  f\mu\right)  \right\}  \text{ .}%
\]
Note that $\int_{\Omega}u(w)f\,d\mu\in(-\infty,\infty]$ for each
$w\in\mathcal{A}\left(  f\right)  $. We define
\[
H_{u}\left(  f\right)  :=\ln u^{-1}\left(  N_{u}\left(  f\right)  \right)
\]
and call it the $u$-\textbf{entropy} (\textbf{utility maximising entropy)
}of~$f$.
\end{definition}

The next proposition follows immediately from the definitions.

\begin{proposition}
\label{ordrel}Let $\mu\in M_{1}$ and $f\in D(\mu)$. Then
\begin{align*}
N_{u}\left(  f\right)   &  =N_{u}\left(  f\mu\parallel\mu\right)  \text{ ,}\\
H_{u}\left(  f\right)   &  =H_{u}\left(  f\mu\parallel\mu\right)  \text{ .}%
\end{align*}
\end{proposition}

\subsubsection{Properties}

\begin{proposition}
\label{entbasbou}The following inequalities hold:
\begin{align*}
u\left(  1\right)   &  \leq N_{u}\left(  \nu\parallel\mu\right)  \leq u\left(
\infty\right)  \text{ ,}\\
0  &  \leq H_{u}\left(  \nu\parallel\mu\right)  \leq\infty\text{ .}%
\end{align*}
\end{proposition}

\begin{proof}
Taking $w\equiv1\in\mathcal{A}\left(  \nu,\mu\right)  ,$ we obtain the lower
bound. The upper bound follows immediately from the definition.
\end{proof}

\begin{proposition}
\label{conentnum}Let $\mu,\nu_{1},\nu_{2}\in M_{1}$ and $a\in\left[
0,1\right]  $. Then
\[
N_{u}\left(  a\nu_{1}+\left(  1-a\right)  \nu_{2}\parallel\mu\right)  \leq
aN_{u}\left(  \nu_{1}\parallel\mu\right)  +\left(  1-a\right)  N_{u}\left(
\nu_{2}\parallel\mu\right)  \text{ .}%
\]
\end{proposition}

\begin{proof}
Put $\nu:=a\nu_{1}+\left(  1-a\right)  \nu_{2}$. First observe that for $w\in
D\left(  \mu\right)  $ we have $\int_{\Omega}u^{-}\left(  w\right)  d\nu
=a\int_{\Omega}u^{-}\left(  w\right)  d\nu_{1}+\left(  1-a\right)
\int_{\Omega}u^{-}\left(  w\right)  d\nu_{2}$, and so $\mathcal{A}\left(
\nu,\mu\right)  =\mathcal{A}\left(  \nu_{1},\mu\right)  \cap\mathcal{A}\left(
\nu_{2},\mu\right)  $. Hence
\begin{align*}
&  N_{u}\left(  a\nu_{1}+\left(  1-a\right)  \nu_{2}\parallel\mu\right) \\
&  =\sup\left\{  \int_{\Omega}u\left(  w\right)  d\nu:w\in\mathcal{A}\left(
\nu,\mu\right)  \right\} \\
&  =\sup\left\{  a\int_{\Omega}u\left(  w\right)  d\nu_{1}+\left(  1-a\right)
\int_{\Omega}u\left(  w\right)  d\nu_{2}:w\in\mathcal{A}\left(  \nu
,\mu\right)  \right\} \\
&  \leq a\sup\left\{  \int_{\Omega}u\left(  w\right)  d\nu_{1}\in
\mathcal{A}\left(  \nu,\mu\right)  \right\}  +\left(  1-a\right)  \sup\left\{
\int_{\Omega}u\left(  w\right)  d\nu_{2}:w\in\mathcal{A}\left(  \nu
,\mu\right)  \right\} \\
&  \leq a\sup\left\{  \int_{\Omega}u\left(  w\right)  d\nu_{1}\in
\mathcal{A}\left(  \nu_{1},\mu\right)  \right\}  +\left(  1-a\right)
\sup\left\{  \int_{\Omega}u\left(  w\right)  d\nu_{2}:w\in\mathcal{A}\left(
\nu_{2},\mu\right)  \right\} \\
&  =aN_{u}\left(  \nu_{1}\parallel\mu\right)  +\left(  1-a\right)
N_{u}\left(  \nu_{2}\parallel\mu\right)  \text{ ,}%
\end{align*}
as desired.
\end{proof}

Next we show that relative $u$-entropy can be reduced to the case when $\nu
\ll\mu$.

\begin{theorem}
\label{relord}Let $\mu,\nu\in M_{1}$. Then
\[
N_{u}\left(  \nu\parallel\mu\right)  =\nu_{\perp}\left(  \Omega\right)
u\left(  \infty\right)  +\nu_{\ll}\left(  \Omega\right)  N_{u}\left(
\frac{\nu_{\ll}}{\nu_{\ll}\left(  \Omega\right)  }\parallel\mu\right)  \text{
,}%
\]
where $\nu_{\perp}+\nu_{\ll}=\nu$ is the Lebesgue decomposition of $\nu$ into
the singular part $\nu_{\perp}$ and absolutely continuous part $\nu_{\ll}$
with respect to~$\mu$.
\end{theorem}

\begin{proof}
Let $A\in\Sigma$ be such that $\mu\left(  A\right)  =0$ and $\nu_{\perp
}\left(  A\right)  =\nu_{\perp}\left(  \Omega\right)  $.

\emph{Step 1.} If $\nu\ll\mu$, that is, $\nu_{\perp}\left(  \Omega\right)
=0$, then the assertion is trivial. Suppose that $\nu\perp\mu$, i.e.,
$\nu_{\perp}\left(  \Omega\right)  =1$. Then $\nu\left(  A\right)  =1$ and
$w_{n}:\equiv n1_{A}+1_{A^{\urcorner}}\in\mathcal{A}\left(  \nu,\mu\right)  $
for $n\in\mathbb{N}$, and $\int_{\Omega}u\left(  w_{n}\right)  d\nu=u\left(
n\right)  $. Hence and from Proposition \ref{entbasbou}$\ $we get
$N_{u}\left(  \nu\parallel\mu\right)  =u\left(  \infty\right)  $, as required.

\emph{Step 2.} Now we assume that $0<\nu_{\perp}\left(  \Omega\right)  <1$.
Note that $\nu=\nu_{\perp}\left(  \Omega\right)  \frac{\nu_{\perp}}{\nu
_{\perp}\left(  \Omega\right)  }+\nu_{\ll}\left(  \Omega\right)  \frac
{\nu_{\ll}}{\nu_{\ll}\left(  \Omega\right)  }$, and from
Proposition~\ref{conentnum} and from Step~1 we get
\[
N_{u}\left(  \nu\parallel\mu\right)  \leq\nu_{\perp}\left(  \Omega\right)
u\left(  \infty\right)  +\nu_{\ll}\left(  \Omega\right)  N_{u}\left(
\frac{\nu_{\ll}}{\nu_{\ll}\left(  \Omega\right)  }\parallel\mu\right)  \text{
.}%
\]
Let now $w\in\mathcal{A}\left(  \frac{\nu_{\ll}}{\nu_{\ll}\left(
\Omega\right)  },\mu\right)  $. Put $w_{n}:\equiv n1_{A}+w1_{A^{\urcorner}}$
for $n\in\mathbb{N}$. Clearly, $w_{n}\in\mathcal{A}\left(  \nu,\mu\right)  $
and
\[
\nu_{\perp}\left(  \Omega\right)  u\left(  n\right)  +\nu_{\ll}\left(
\Omega\right)  \int_{\Omega}u\left(  w\right)  d\frac{\nu_{\ll}}{\nu_{\ll
}\left(  \Omega\right)  }=\int_{\Omega}u\left(  w_{n}\right)  d\nu\leq
N_{u}\left(  \nu\parallel\mu\right)  \text{ .}%
\]
Taking $n\rightarrow\infty$ completes the proof.
\end{proof}

\begin{corollary}
In particular, if $u\left(  \infty\right)  =\infty$ and $\nu$ is not
absolutely continuous with respect to~$\mu$, or if $\nu\perp\mu$, then
$N_{u}\left(  \nu\parallel\mu\right)  =u\left(  \infty\right)  $ and
$H_{u}\left(  \nu\parallel\mu\right)  =\infty$.
\end{corollary}

\begin{proposition}
\label{utiinffin}Let $\mu,\nu\in M_{1}$. Then the following conditions are equivalent:

\begin{enumerate}
\item [$(1)$]$N_{u}\left(  \nu\parallel\mu\right)  <u\left(  \infty\right)  $ ;

\item[$(2)$] $N_{u}\left(  \nu\parallel\mu\right)  <\infty$ ;

\item[$(3)$] $H_{u}\left(  \nu\parallel\mu\right)  <\infty$ .
\end{enumerate}

\noindent In particular, all three conditions are satisfied for any utility
function $u$ such that $u\left(  \infty\right)  <\infty$.
\end{proposition}

\begin{proof}
The implications $(1)\Rightarrow(3)\Rightarrow(2)$ are obvious, as is
$(2)\Rightarrow(1)$ when $u(\infty)=\infty$.

We shall prove $2)\Rightarrow1)$ when $u(\infty)<\infty$. Put $A_{n}:=\left\{
w\geq n\right\}  $ for any $n\in\mathbb{N}$. Then $\bigcap_{n\in\mathbb{N}%
}A_{n}=\emptyset$. Consequently, there exists an $n\in\mathbb{N}$ such that
$\nu\left(  A_{n}\right)  =:\gamma<1$. Hence
\begin{align*}
\int_{\Omega}u\left(  w\right)  d\nu &  =\int_{A_{n}}u\left(  w\right)
d\nu+\int_{\left(  A_{n}\right)  ^{\urcorner}}u\left(  w\right)  d\nu\\
&  \leq u\left(  \infty\right)  \gamma+u\left(  n\right)  \left(
1-\gamma\right)
\end{align*}
for any $w\in\mathcal{A}\left(  \nu,\mu\right)  $. Thus $N_{u}\left(
\nu\parallel\mu\right)  \leq u\left(  \infty\right)  \gamma+u\left(  n\right)
\left(  1-\gamma\right)  <u\left(  \infty\right)  $, as required.
\end{proof}

\begin{proposition}
The following conditions are equivalent:

\begin{enumerate}
\item [$(1)$]$H_{u}\left(  \nu\parallel\mu\right)  =0$ ;

\item[$(2)$] $\nu=\mu$ .
\end{enumerate}
\end{proposition}

\begin{proof}
$(1)\Rightarrow(2).$ Let $\nu=\mu$. Take $w\in\mathcal{A}\left(  \nu
,\mu\right)  $. By Jensen's inequality $\int_{\Omega}u(w)d\nu\leq u\left(
\int_{\Omega}wd\nu\right)  =u\left(  1\right)  $. Hence, by
Proposition~\ref{entbasbou}, $N_{u}\left(  \nu\parallel\mu\right)  =u(1)$, and
so $H_{u}\left(  \nu\parallel\mu\right)  =0$.

$(2)\Rightarrow(1).$ Suppose that $\nu\neq\mu$. Then there is an $A\in\Sigma$
such that $\mu(A)\neq\nu(A)$ and $\mu(A)<1$. We put
\begin{align*}
w_{a}  &  :=a1_{A}+\frac{1-a\mu(A)}{\mu(A^{\urcorner})}1_{A^{\urcorner}%
}~\text{,}\\
\varphi\left(  a\right)   &  :=\int_{\Omega}u(w_{a})d\nu=u(a)\nu(A)+u\left(
\frac{1-a\mu(A)}{\mu(A^{\urcorner})}\right)  \nu(A^{\urcorner})
\end{align*}
for any $a\in\left(  0,1/\mu(A)\right)  $. Clearly, $w_{a}\in\mathcal{A}%
\left(  \nu,\mu\right)  $ and $w_{1}\equiv1$. Moreover, $\varphi^{\prime
}\left(  1\right)  =u^{\prime}\left(  1\right)  \frac{\nu(A)-\mu(A)}%
{\mu(A^{\urcorner})}\neq0$. Hence there exists an $a\in\left(  0,1/\mu
(A)\right)  $ such that $\int_{\Omega}u\left(  w_{a}\right)  d\nu
=\varphi\left(  a\right)  >\varphi\left(  1\right)  =\int_{\Omega}u\left(
w_{1}\right)  d\nu=u\left(  1\right)  $. Thus $H_{u}\left(  \nu\parallel
\mu\right)  >0$.
\end{proof}

\begin{proposition}
[linear transformation]\label{lintraent}Let $u:(0,\infty)\mathbb{\rightarrow
}\mathbb{R}$ be a utility function, let $a>0$ and let $b\in\mathbb{R}$. Then
$\widetilde{u}=au+b$ is a utility function, and for any $\nu,\mu\in M_{1}$%
\begin{align*}
N_{\widetilde{u}}\left(  \nu\parallel\mu\right)   &  =aN_{u}\left(
\nu\parallel\mu\right)  +b\text{ ,}\\
H_{\widetilde{u}}\left(  \nu\parallel\mu\right)   &  =H_{u}\left(
\nu\parallel\mu\right)  \text{ .}%
\end{align*}
\end{proposition}

\begin{proof}
This follows immediately from the definition.
\end{proof}

\begin{remark}
\upshape It has recently been proved by Urba\'{n}ski \cite{Urb07} that in
probability spaces without atoms $\widetilde{u}=au+b$ is not only a sufficient
condition, but in fact an equivalent condition for $H_{\widetilde{u}}=H_{u}$.
The equivalence can fail in a probability space with atoms.
\end{remark}

In \cite[Theorem 20]{SloZas01} we established a formula for $u$-entropy by
convex duality methods. Namely, under the reasonable asymptotic elasticity
assumption, if $f\in D\left(  \mu\right)  $, then
\begin{align*}
N_{u}\left(  f\right)   &  =\int_{\Omega}u\left(  I\left(  \Lambda
_{f}/f\right)  \right)  f\,d\mu=\int_{\Omega}u^{\ast}\left(  \Lambda
_{f}/f\right)  f\,d\mu+\Lambda_{f}\text{ ,}\\
H_{u}\left(  f\right)   &  =\ln u^{-1}\left(  \int_{\Omega}u\left(  I\left(
\Lambda_{f}/f\right)  \right)  f\,d\mu\right)  \text{ ,}%
\end{align*}
where $\Lambda_{f}>0$ is given implicitly as the unique solution of
\[
\int_{\Omega}I\left(  \Lambda_{f}/f\right)  \,d\mu=1\text{ .}%
\]
Combined with Theorem~\ref{relord}, this makes it possible to evaluate the
relative $u$-entropy $H_{u}\left(  \nu\parallel\mu\right)  $ for any $\nu
,\mu\in M_{1}$.

\begin{example}
[logarithmic utility]\upshape Let $u:(0,\infty)\mathbb{\rightarrow}\mathbb{R}$
be given by $u\left(  x\right)  =\ln x$ for $x\in(0,\infty)$. Then $H_{u}\ $is
equal to the \textbf{Boltzmann-Shannon relative entropy}
\[
H_{1}\left(  \nu\parallel\mu\right)  =\left\{
\begin{tabular}
[c]{ll}%
$\int_{\Omega}\frac{d\nu}{d\mu}\ln\frac{d\nu}{d\mu}d\mu$ & if $\nu\ll\mu$ ,\\
$\infty$ & otherwise
\end{tabular}
\right.
\]
for $\mu,\nu\in M_{1}$.
\end{example}

\begin{example}
[isoelastic utility]\label{exlisouti}\upshape Let $u:(0,\infty
)\mathbb{\rightarrow}\mathbb{R}$ be given by $u\left(  x\right)  =\frac
{1}{\gamma}\left(  x^{\gamma}-1\right)  $ for $\gamma\in\left(  -\infty
,0\right)  \cup\left(  0,1\right)  $ and $x\in(0,\infty)$. Then $H_{u}$ is
equal to the \textbf{R\'{e}nyi relative entropy of order }$\alpha=\left(
1-\gamma\right)  ^{-1}\in\left(  0,1\right)  \cup\left(  1,\infty\right)  $
given by
\[
H_{\alpha}\left(  \nu\parallel\mu\right)  =\left\{
\begin{tabular}
[c]{ll}%
$\frac{1}{\alpha-1}\ln\int_{\Omega}\left(  \frac{d\nu_{\ll}}{d\mu}\right)
^{\alpha}d\mu\smallskip$ & if $\gamma\in\left(  -\infty,0\right)  $
,\smallskip\\
$\frac{1}{\alpha-1}\ln\int_{\Omega}\left(  \frac{d\nu}{d\mu}\right)  ^{\alpha
}d\mu\smallskip$ & if $\gamma\in\left(  0,1\right)  $ and $\nu\ll\mu$
,$\smallskip$\\
$\infty$ & otherwise
\end{tabular}
\right.
\]
for $\mu,\nu\in M_{1}$.
\end{example}

\begin{remark}
\upshape
The Boltzmann-Shannon relative entropy was introduced in\ \cite{KulLei51}
under the name of directed divergence. It is also called the Kullback-Leibler
divergence, relative information, conditional entropy, information gain or
function of discrimination. The definition of the R\'{e}nyi relative entropy
(or divergence) of order $\alpha$ was proposed in \cite{Ren61}.
\end{remark}

\section{Discrete $u$-entropy}

Let $\Omega=\{\omega_{1},\ldots,\omega_{k}\}$ be a finite probability space
equipped with the sigma-field $\Sigma=2^{\Omega}$ of all subsets of~$\Omega$.
The family of probability measures on $(\Omega,\Sigma)$ will be denoted
by~$S_{k}$. For any $p\in S_{k}$ we shall write $p_{i}=p(\omega_{i}) $ for
$i=1,\ldots,k$. Thus, we can identify $S_{k}$ with the set of probability
vectors $\left\{  p\in\mathbb{R}^{k}:\sum_{i=1}^{k}p_{i}=1\text{ and }%
p_{i}\geq0\text{ for }i=1,\ldots,k\right\}  $. Our definition of the relative
$u$-entropy covers also the discrete case. In this situation (though not
necessarily in the general case) it is also possible to define the
(non-relative) $u$-entropy as follows.

\begin{definition}
\label{defdiscrentr}\upshape
Let $u:\left(  0,\infty\right)  \mathbb{\rightarrow}\mathbb{R}$ be a utility
function, and let $p\in S_{k}$. Then we put
\[
n_{u}\left(  p\right)  :=\sup_{w\in S_{k}}\sum_{i=1}^{k}u\left(  w_{i}\right)
p_{i}\text{ ,}%
\]
and define the \textbf{discrete }$u$\textbf{-entropy of }$p$ by
\[
h_{u}\left(  p\right)  :=-\ln u^{-1}\left(  n_{u}\left(  p\right)  \right)
\text{ .}%
\]
\end{definition}

\begin{remark}
\upshape
Note that $h_{u}$ depends only on the restriction of $u$ to $(0,1]$.
\end{remark}

\begin{proposition}
\label{discon}Let $u:(0,\infty)\mathbb{\rightarrow}\mathbb{R}$ be a utility
function. Let $p\in S_{k}$ and let $p_{(k)}\in S_{k}$ be the uniform
probability vector, that is, $(p_{\left(  k\right)  })_{i}=1/k$ for each
$i=1,\ldots,k$. Then
\[
H_{u}\left(  p\parallel p_{(k)}\right)  =\ln k-h_{u_{k}}\left(  p\right)
\text{ ,}%
\]
where $u_{k}:(0,\infty)\mathbb{\rightarrow}\mathbb{R}$ is the rescaled utility
function
\begin{equation}
u_{k}\left(  x\right)  :=u\left(  kx\right)  \label{urescaled}%
\end{equation}
for $x\in(0,\infty)$.
\end{proposition}

\begin{proof}
Since $u^{-1}=ku_{k}^{-1}$ and
\begin{align*}
N_{u}(p\parallel p_{(k)})  &  =\sup_{w\in D(p_{(k)})}\sum_{i=1}^{k}%
u(w_{i})p_{i}\\
&  =\sup_{w\in S_{k}}\sum_{i=1}^{k}u(kw_{i})p_{i}=\sup_{w\in S_{k}}\sum
_{i=1}^{k}u_{k}(w_{i})p_{i}=n_{u_{k}}(p)\text{ ,}%
\end{align*}
it follows that
\[
H_{u}\left(  p\parallel p_{(k)}\right)  =\ln u^{-1}(N_{u}(p\parallel
p_{(k)}))=\ln[ku_{k}^{-1}(n_{u_{k}}(p))]=\ln k-h_{u_{k}}(p)\text{ .}%
\]
\end{proof}

Using the above statement we can deduce many properties of discrete
$u$-entropy from the respective properties of relative $u$-entropy. However,
one can also prove them straightforwardly without assuming anything about the
behaviour of the function $u$ outside the interval $\left(  0,1\right]  $. We
could assume that $u:\left(  0,1\right]  \mathbb{\rightarrow}\mathbb{R}$ is a
strictly concave strictly increasing continuously\ differentiable function
such that $\lim_{x\searrow0}u^{\prime}\left(  x\right)  =\infty$. In this case
$I:=\left(  u^{\prime}\right)  ^{-1}$ would be defined on the interval
$\left[  u^{\prime}\left(  1\right)  ,\infty\right)  $. The proofs of the
following properties of discrete $u$-entropy are elementary.

\begin{proposition}
Let $p=\left(  p_{1},\ldots,p_{k}\right)  \in S_{k}$. Then

\begin{enumerate}
\item [$(1)$]$0\leq h_{u}(p)\leq\ln k$.

\item[$(2)$] $h_{u}(p)=0$ iff $p_{i}=1$ for some $i=1,\ldots,k$.

\item[$(3)$] $h_{u}(p)=\ln k$ iff $p=p_{\left(  k\right)  }$.

\item[$(4)$] $h_{u}(p)=h_{u}((p_{\pi\left(  1\right)  },\ldots,p_{\pi\left(
k\right)  }))$ for every permutation $\pi$.

\item[$(5)$] $h_{u}(p)=h_{u}(\left(  p_{1},\ldots,p_{k},0\right)  )$.

\item[$(6)$] For $a>0$ and $b\in\mathbb{R}$ we have $h_{au+b}=h_{u}$.
\end{enumerate}
\end{proposition}

The proof of the formula for $u$-entropy in the discrete case is also
elementary and, by contrast to the general case, it does not require any
further assumptions on $u$.

\begin{proposition}
\label{disentfor}Let $p\in S_{k}$. Then:

\begin{enumerate}
\item [$(1)$]There exists a unique $\Lambda_{p}\geq\alpha\left(  p\right)
:=u^{\prime}\left(  1\right)  \max\limits_{j=1,\ldots,k}p_{j}>0$ such that
\[
\sum_{i=1}^{k}I\left(  \Lambda_{p}/p_{i}\right)  =1\text{ .}%
\]

\item[$(2)$] The following formulae hold:
\begin{align}
n_{u}\left(  p\right)   &  =\sum_{i=1}^{k}u\left(  I\left(  \Lambda_{p}%
/p_{i}\right)  \right)  p_{i}=\sum_{i=1}^{k}u^{\ast}\left(  \Lambda_{p}%
/p_{i}\right)  p_{i}+\Lambda_{p}\text{ ,}\label{disforent}\\
h_{u}\left(  p\right)   &  =\ln u^{-1}\left(  \sum_{i=1}^{k}u\left(  I\left(
\Lambda_{p}/p_{i}\right)  \right)  p_{i}\right)  \text{ .}\nonumber
\end{align}
\end{enumerate}
\end{proposition}

\begin{proof}
To prove (1) consider the function $\phi_{p}:\left[  \alpha\left(  p\right)
,\infty\right)  \rightarrow(0,\infty)$ given by $\phi_{p}\left(
\Lambda\right)  =\sum_{i=1}^{k}I\left(  \Lambda/p_{i}\right)  $ for
$\Lambda\geq\alpha\left(  p\right)  $. Clearly, $\phi_{p}$ is continuous,
strictly decreasing and satisfies $\phi_{p}\left(  \alpha\left(  p\right)
\right)  =\sum_{i=1}^{k}I\left(  u^{\prime}\left(  1\right)  \max
_{j=1,\ldots,k}p_{j}/p_{i}\right)  \geq1$ and $\lim_{\Lambda\rightarrow\infty
}\phi_{p}\left(  \Lambda\right)  =0$. As a result, there is a unique
$\Lambda_{p}\geq\alpha\left(  p\right)  $ such that $\phi_{p}\left(
\Lambda_{p}\right)  =1$, as required. It follows from (1) that $n_{u}\left(
p\right)  \geq\sum_{i=1}^{k}u\left(  I\left(  \Lambda_{p}/p_{i}\right)
\right)  p_{i}$. To prove the reverse inequality take $w\in S_{k}$. Let
$i=1,\ldots,k$. From the well-known formula $u^{\ast}(y)=u(I(y))-yI(y)$ for
the convex dual we get
\begin{equation}
u\left(  w_{i}\right)  -\left(  \Lambda_{p}/p_{i}\right)  w_{i}\leq u^{\ast
}\left(  \left(  \Lambda_{p}/p_{i}\right)  \right)  =u\left(  I\left(  \left(
\Lambda_{p}/p_{i}\right)  \right)  \right)  -\left(  \Lambda_{p}/p_{i}\right)
I\left(  \left(  \Lambda_{p}/p_{i}\right)  \right)  \text{ .} \label{ine1}%
\end{equation}
Multiplying (\ref{ine1}) by $p_{i}$, summing over $i=1,\ldots,k$, and
adding~$\Lambda_{p}$, we obtain
\[
\sum_{i=1}^{k}u\left(  w_{i}\right)  p_{i}\leq\sum_{i=1}^{k}u^{\ast}\left(
\left(  \Lambda_{p}/p_{i}\right)  \right)  p_{i}+\Lambda_{p}=\sum_{i=1}%
^{k}u\left(  I\left(  \left(  \Lambda_{p}/p_{i}\right)  \right)  \right)
p_{i}%
\]
Taking the supremum of the left-hand side over all such $w$'s, we obtain the assertion.
\end{proof}

\begin{proposition}
The function $h_{u}:S_{k}\rightarrow\left[  0,\ln k\right]  $ is continuous.
\end{proposition}

\begin{proof}
According to (\ref{disforent}), it is enough to prove that $S_{k}\ni
p\rightarrow\Lambda_{p}\in\left[  \alpha\left(  p\right)  ,\infty\right)  $ is
continuous. Define $F:\left\{  \left(  p,\Lambda\right)  \right\}  :p\in
S_{k},\Lambda\in\left[  \alpha\left(  p\right)  ,\infty\right)  \}\rightarrow
\mathbb{R}$ by $F\left(  p,\Lambda\right)  =\sum_{i=1}^{k}I\left(
\Lambda/p_{i}\right)  -1$. Clearly, $F$ is continuous, $F\left(  p,\Lambda
_{p}\right)  =0$ for $p\in S_{k}$, and $\left[  \alpha\left(  p\right)
,\infty\right)  \ni\Lambda\rightarrow F\left(  p,\Lambda\right)  \in
\mathbb{R}$ is strictly decreasing for each $p\in S_{k}$. Now, the assertion
follows from the implicit function theorem for continuous functions.
\end{proof}

\begin{example}
[logarithmic utility]\upshape
Let $u:(0,\infty)\mathbb{\rightarrow}\mathbb{R}$ be given by $u\left(
x\right)  =\ln x$ for $x\in(0,\infty)$. Then the relative $u$-entropy
$H_{u}\left(  p\parallel q\right)  $ is equal to the \textbf{discrete
Boltzmann-Shannon relative entropy (Kullback-Leibler divergence)}
\[
H_{1}\left(  p\parallel q\right)  =\left\{
\begin{tabular}
[c]{ll}%
$\sum\limits_{\substack{i=1,\ldots,k \\q_{i}\neq0}}p_{i}\ln\frac{p_{i}}{q_{i}%
}$ & if $p\ll q$ ,\\
$\infty$ & otherwise
\end{tabular}
\right.
\]
for $p,q\in S_{k}$, and the discrete $u$-entropy $h_{u}(p)$ is equal to the
\textbf{discrete Boltzmann-Shannon entropy}
\[
h_{1}\left(  p\right)  =-\sum\limits_{i=1}^{k}p_{i}\ln p_{i}%
\]
for $p\in S_{k}$.
\end{example}

\begin{example}
[isoelastic utility]\upshape
Let $u:(0,\infty)\mathbb{\rightarrow}\mathbb{R}$ be given by $u\left(
x\right)  =\frac{1}{\gamma}\left(  x^{\gamma}-1\right)  $ for $\gamma
\in\left(  -\infty,0\right)  \cup\left(  0,1\right)  $ and $x\in(0,\infty)$.
Then the discrete relative $u$-entropy $H_{u}\left(  p\parallel p_{(k)}%
\right)  $ is equal to the \textbf{discrete R\'{e}nyi relative entropy
(divergence) of order }$\alpha=\left(  1-\gamma\right)  ^{-1}\in\left(
0,1\right)  \cup\left(  1,\infty\right)  $%
\[
h_{\alpha}\left(  p\parallel q\right)  =\left\{
\begin{tabular}
[c]{ll}%
$\frac{1}{\alpha-1}\ln\sum\limits_{i=1,\ldots,k}p_{i}^{\alpha}q_{i}^{1-\alpha
}$ & if $\gamma\in\left(  -\infty,0\right)  $ ,\smallskip\\
$\frac{1}{\alpha-1}\ln\sum\limits_{\substack{i=1,\ldots,k\\q_{i}\neq0}%
}p_{i}^{\alpha}q_{i}^{1-\alpha}$ & if $\gamma\in\left(  0,1\right)  $ and
$p\ll q$ ,\\
$\infty$ & otherwise
\end{tabular}
\right.
\]
for $p,q\in S_{k}$, and the discrete $u$-entropy $h_{u}(p)$ is equal to the
\textbf{discrete R\'{e}nyi entropy of order }$\alpha$%
\[
h_{\alpha}\left(  p\right)  =\frac{1}{1-\alpha}\ln\sum\limits_{i=1}^{k}%
p_{i}^{\alpha}%
\]
for $p\in S_{k}$.
\end{example}

\section{Relationships to other utility based concepts of entropy}

\subsection{Friedman-Huang-Sandow $U$\textit{-}entropy}

In~\cite{FriHuaSan07} (see also~\cite{FriHuaSan05}) the authors defined two
quantities, which they called the $U$\textbf{-entropy}\textsl{\ }and
$U$\textbf{-relative entropy}, noting their similarity to the $u$-entropy
defined (in a much more general setting) in~\cite{SloZas04}. In fact the
$U$-entropy and $U$-relative entropy of Friedman, Huang and
Sandow~\cite{FriHuaSan07}, \cite{FriHuaSan05} can be reduced by a simple
transformation to the relative $u$-entropy discussed in the present paper, and
so to the $u$-entropy defined in~\cite{SloZas04}. As a result, the properties
of $U$-entropy and $U$-relative entropy claimed in~\cite{FriHuaSan05},
\cite{FriHuaSan07} turn out to be immediate corollaries of the results
of~\cite{SloZas04}, as shown below.

In the notation of the present paper the definitions in~\cite{FriHuaSan07}
take the following form.

\begin{definition}
\upshape
(Definition~5 from~\cite{FriHuaSan07}) Let $u:(0,\infty)\mathbb{\rightarrow
}\mathbb{R}$ be a utility function and let $p,q\in S_{k}$. If $p\ll q$, then
the Friedman-Huang-Sandow $U$\textbf{-relative entropy of }$p$\textbf{\ with
respect to~}$q$ is defined by
\[
D_{u}\left(  p\parallel q\right)  :=\sup_{w\in S_{k}}\sum_{i=1}^{k}u\left(
\frac{w_{i}}{q_{i}}\right)  p_{i}-u\left(  1\right)  \text{ .}%
\]
\end{definition}

\begin{remark}
\upshape
By contrast to \cite{FriHuaSan07}, it is not assumed here that $u(1)=0.$ To
compensate, we subtract $u(1)$ on the right-hand side of the formula defining
$D_{u}\left(  p\parallel q\right)  $. The same applies to the formula defining
$H_{u}(p)$ below. Moreover, instead of $w\in S_{k}$ it is only assumed
in~\cite{FriHuaSan07} that $\sum_{i=1}^{k}w_{i}=1$, but presumably there is
also a silent assumption that $w_{i}/q_{i}$ belongs to the domain of~$u$ for
each~$i$. In our case this means that, additionally, $w\geq0$, so that $w\in
S_{k}$. The definitions and results easily extend to utility functions defined
on an interval $(a,b)$ other than $(0,\infty)$. If $p$ is not absolutely
continuous with respect to~$q$, then $D_{u}\left(  p\parallel q\right)  $ is undefined.
\end{remark}

\begin{remark}
\upshape
The relative entropy defined in~\cite{FriHuaSan07} coincides with the
\textbf{decision maker's optimal expected utility}\textsl{\ }introduced in
\cite[p.13]{JosNauWin07}.
\end{remark}

\begin{definition}
\label{Frieeralentropydef}\upshape
(Definition 6 from \cite{FriHuaSan07}) Let $u:(0,\infty)\mathbb{\rightarrow
}\mathbb{R}$ be a utility function and let $p\in S_{k}$. Then the
Friedman-Huang-Sandow $U$\textbf{-entropy of }$p$ is defined by
\[
H_{u}\left(  p\right)  =u\left(  k\right)  -u\left(  1\right)  -D_{u}\left(
p\parallel p_{(k)}\right)  \text{ .}%
\]
\end{definition}

\begin{proposition}
\label{SloZasFri}Let $u:(0,\infty)\mathbb{\rightarrow}\mathbb{R}$ be a utility
function. Then the following properties hold:

\begin{enumerate}
\item [$(1)$]For any $p,q\in S_{k}$ such that $p\ll q$%
\[
D_{u}\left(  p\parallel q\right)  =u\left(  e^{H_{u}\left(  p\parallel
q\right)  }\right)  -u\left(  1\right)  \text{ .}%
\]

\item[$(2)$] For any $p\in S_{k}$%
\[
H_{u}\left(  p\right)  =u_{k}\left(  1\right)  -u_{k}\left(  e^{-h_{u_{k}%
}\left(  p\right)  }\right)  \text{ ,}%
\]
where $u_{k}$ is the rescaled utility function defined by~$(\ref{urescaled})$.
\end{enumerate}
\end{proposition}

\begin{proof}
$(1)$ If $p\ll q$, then $wq\in S_{k}$ is equivalent to $w\in\mathcal{A}(p,q)$.
Hence
\begin{align*}
D_{u}\left(  p\parallel q\right)   &  =\sup_{w\in S_{k}}\sum_{i=1}^{k}u\left(
w_{i}/q_{i}\right)  p_{i}-u\left(  1\right) \\
&  =\sup_{w\in\mathcal{A}(p,q)}\sum_{i=1}^{k}u\left(  w_{i}\right)
p_{i}-u\left(  1\right)  =N_{u}\left(  p\parallel q\right)  -u\left(
1\right)  \text{ .}%
\end{align*}
The claim follows since $N_{u}\left(  p\parallel q\right)  =u\left(
e^{H_{u}\left(  p\parallel q\right)  }\right)  $.

$(2)$ This follows immediately from~$(1)$ and Proposition~\ref{discon}:
\begin{align*}
H_{u}\left(  p\right)   &  =u\left(  k\right)  -u\left(  1\right)
-D_{u}\left(  p\parallel p_{(k)}\right)  =u\left(  k\right)  -u\left(
e^{H_{u}\left(  p\parallel p_{(k)}\right)  }\right) \\
&  =u\left(  k\right)  -u\left(  e^{\ln k-h_{u_{k}}\left(  p\right)  }\right)
=u_{k}\left(  1\right)  -u_{k}\left(  e^{-h_{u_{k}}\left(  p\right)  }\right)
\text{ .}%
\end{align*}
\end{proof}

The $U$-relative entropy $D_{u}\left(  p\parallel q\right)  $ and the
$U$-entropy $H_{u}\left(  p\right)  $ of Friedman Huang and Sandow
\cite{FriHuaSan07} are therefore related to the $u$-entropy of~\cite{SloZas04}
by
\[
D_{u}\left(  p\parallel q\right)  =u\left(  e^{H_{u}\left(  \frac{dp}%
{dq}\right)  }\right)  -u\left(  1\right)
\]
for each $p,q\in S_{k}$ such that $p\ll q$, and by
\[
H_{u}\left(  p\right)  =u\left(  k\right)  -u\left(  e^{H_{u}\left(  \frac
{dp}{dp_{(k)}}\right)  }\right)
\]
for each $p\in S_{k}$. Because of this, the following results
in~\cite{FriHuaSan07} are immediate consequences of the corresponding earlier
results in~\cite{SloZas04}:
\[%
\begin{tabular}
[c]{lll}%
\cite{SloZas04}\smallskip &  & \cite{FriHuaSan07}\smallskip\\
Theorem 20 & $\Rightarrow$ & Lemma 1\\
Propositions 8 and 10 & $\Rightarrow$ & Corollary 1.(i)\\
Proposition 13 & $\Rightarrow$ & Corollary 1.(ii), (iv)\\
Theorem 23 & $\Rightarrow$ & Corollaries 3, 6 and 7
\end{tabular}
\]

\begin{remark}
\upshape The results of \cite{SloZas04} are valid in a much more general
situation of arbitrary probability spaces, which requires the asymptotic
elasticity assumption
\[
\mathrm{AE}(u):=\limsup_{x\nearrow\infty}\frac{xu^{\prime}(x)}{u(x)}<1
\]
to hold. In the discrete case this assumption is unnecessary and all the
arguments in~\cite{SloZas04} work without it.
\end{remark}

\subsection{Arimoto entropy}

A similar construction of entropy was first proposed by Arimoto~\cite{Ari71}
without any explicit reference to the notion of utility.

\begin{definition}
\upshape(Arimoto~\cite{Ari71}, see also~\cite{Tan_s39}) For a non-negative
function $f:(0,1]\rightarrow\mathbb{R}$ continuously differentiable on $(0,1]
$ and such that $f(1)=0$ Arimoto's entropy is defined by
\[
H_{f}^{A}\left(  p\right)  :=\inf_{w\in S_{k}}\sum_{i=1}^{k}f\left(
w_{i}\right)  p_{i}%
\]
for $p\in S_{k}$.
\end{definition}

This was further generalized in~\cite{ShaSon74} and also interpreted in~
\cite[Example 6]{MorParVaj96} in terms of \textsl{prior Bayes risk}, where $f
$ plays the role of an individual uncertainty function. Arimoto's entropy is
related to the entropy~$H_{u}$ defined in~\cite{FriHuaSan07} (see
Definition~\ref{Frieeralentropydef} above) and to $h_{u}$
(Definition~\ref{defdiscrentr}) as follows.

\begin{proposition}
\label{Ari}Let $u:(0,\infty)\mathbb{\rightarrow}\mathbb{R}$ be a utility
function such that $u\left(  1\right)  =0$, and let $p\in S_{k}$. Then
\[
H_{-u}^{A}\left(  p\right)  =H_{u_{1/k}}\left(  p\right)  =-u\!\left(
e^{-h_{u}\left(  p\right)  }\right)  .
\]
\end{proposition}

\begin{proof}
This follows immediately from the definitions and Proposition~\ref{SloZasFri}%
~$(2)$:
\begin{align*}
H_{-u}^{A}\left(  p\right)   &  =\inf_{w\in S_{k}}\sum_{i=1}^{k}\left[
-u\left(  w_{i}\right)  \right]  p_{i}=-\sup_{w\in S_{k}}\sum_{i=1}%
^{k}u\left(  w_{i}\right)  p_{i}\\
&  =-n_{u}(p)=-u\!\left(  e^{-h_{u}\left(  p\right)  }\right)  =H_{u_{1/k}%
}\left(  p\right)  \text{ .}%
\end{align*}
\end{proof}

\begin{example}
[logarithmic utility]\upshape Let $u$ be the logarithmic utility. For $p,q\in
S_{k}$ and $p\ll q$ we have $D_{u}\left(  p\parallel q\right)  =h_{_{1}%
}\left(  p\parallel q\right)  $ and $H_{u}\left(  p\right)  =H_{-u}^{A}\left(
p\right)  =h_{_{1}}\left(  p\right)  $.
\end{example}

\begin{example}
[isoelastic utility]\upshape Let $u$ be the isoelastic utility of order
$\gamma\in\left(  -\infty,0\right)  \cup\left(  0,1\right)  $ and let
$\alpha=\left(  1-\gamma\right)  ^{-1}$.

\begin{enumerate}
\item  For $p,q\in S_{k}$ with $p\ll q$%
\[
D_{u}\left(  p\parallel q\right)  =\frac{\alpha}{\alpha-1}\left(  \left(
\sum\limits_{i=1,\ldots,k}p_{i}^{\alpha}q_{i}^{1-\alpha}\right)  ^{\frac
{1}{\alpha}}-1\right)
\]
is proportional to the \textbf{Sharma-Mittal relative entropy of order
}$\alpha$\textbf{\ and degree }$2-1/a$ ;

\item  For $p\in S_{k}$%
\begin{align*}
H_{u}\left(  p\right)   &  =k^{\frac{\alpha-1}{\alpha}}\frac{\alpha}{1-\alpha
}\left(  \left(  \sum\limits_{i=1}^{k}p_{i}^{\alpha}\right)  ^{\frac{1}%
{\alpha}}-1\right) \\
&  =k^{\frac{\alpha-1}{\alpha}}H_{-u}^{A}\text{ ,}%
\end{align*}
where $H_{-u}^{A}$ is called the \textbf{Arimoto entropy of kind }$1/\alpha$.
\end{enumerate}
\end{example}

\begin{remark}
\upshape The Sharma-Mittal relative entropy was introduced in~\cite{ShaMit75}.
In~\cite{JosNauWin07} the Sharma-Mittal relative entropy of order~$\alpha$ and
degree $2-1/a$ is called the \textbf{pseudospherical divergence of
order}~$\alpha$. The Arimoto entropy was introduced by Arimoto~\cite{Ari71}
and further elaborated in~\cite{BoeLub80}, see also~\cite{Tan_s39}.
\end{remark}

\subsection{Frittelli generalised distance}

This notion of generalised distance in the set of probability measures was
introduced in~\cite{Fri00a} as a tool for solving the convex dual problem to
that of computing the value of a financial security consistent with the
no-arbitrage principle in an incomplete market model in a utility maximisation framework.

\begin{definition}
\upshape
(Definition~9 and formula~(9) from~\cite{Fri00a}) Let $\mu,\nu\in M^{1}$ be
such that $\mu\ll\nu$. Then put
\[
\Delta_{u}(\mu,\nu):=\sup_{\Lambda>0}\left(  \Lambda+\int_{\Omega}u^{\ast
}\left(  \Lambda\frac{d\mu}{d\nu}\right)  d\nu\right)  \text{ ,}%
\]
where $u^{\ast}$ is the convex dual to the utility function~$u$ given
by~(\ref{defconuti}) and (\ref{defconinf}), and define the \textbf{Frittelli
generalised distance} by
\[
\delta_{u}(\mu,\nu):=u^{-1}\left(  \Delta_{u}(\mu,\nu)\right)  -1\text{ .}%
\]
\end{definition}

\begin{remark}
\upshape
Note the different (but equivalent) conventions as compared to~\cite{Fri00a}.
The differences lie in using the supremum rather than infimum coupled with
different signs of certain expressions in the definitions of $u^{\ast}$ and
$\Delta_{u}(\mu,\nu)$. The quantities $\Delta_{u}(\mu,\nu)$ and $\delta
_{u}(\mu,\nu)$ are denoted by $\Delta(\mu,\nu;1)$ and $\delta(\mu,\nu;1)$ in~
\cite{Fri00a}.
\end{remark}

\begin{proposition}
\label{Deltadecomp}For any $\mu,\nu\in M^{1}$ such that $\mu\ll\nu$%
\[
\Delta_{u}(\mu,\nu)=u(\infty)\nu_{\perp}(\Omega)+\nu_{\ll}(\Omega)\Delta
_{u}\!\left(  \mu,\frac{\nu_{\ll}}{\nu_{\ll}(\Omega)}\right)  \text{ ,}%
\]
where $\nu_{\perp}+\nu_{\ll}=\nu$ is the Lebesgue decomposition of $\nu$ into
the singular part $\nu_{\perp}$ and absolutely continuous part $\nu_{\ll}$
with respect to~$\mu$.
\end{proposition}

\begin{proof}
Because $\mu\ll\nu$, it follows that $\frac{d\mu}{d\nu}=0$ a.s.\ with respect
to~$\nu_{\perp}$ and $\nu_{\ll}(\Omega)>0$. We can assume that $\nu_{\perp
}(\Omega)>0$, since otherwise the assertion is obvious. Put $\tilde{\nu}_{\ll
}=\frac{\nu_{\ll}}{\nu_{\ll}(\Omega)}$. As a result,
\begin{align*}
\Delta_{u}(\mu,\nu)  &  =\sup_{\Lambda>0}\left(  \Lambda+\int_{\Omega}u^{\ast
}\left(  \Lambda\frac{d\mu}{d\nu}\right)  d\nu\right) \\
&  =\sup_{\Lambda>0}\left(  \Lambda+\int_{\Omega}u^{\ast}\left(  \Lambda
\frac{d\mu}{d\nu}\right)  d\nu_{\perp}+\int_{\Omega}u^{\ast}\left(
\Lambda\frac{d\mu}{d\nu}\right)  d\nu_{\ll}\right) \\
&  =\nu_{\perp}(\Omega)u^{\ast}\left(  0\right)  +\nu_{\ll}(\Omega
)\sup_{\Lambda>0}\left(  \frac{\Lambda}{\nu_{\ll}(\Omega)}+\int_{\Omega
}u^{\ast}\left(  \frac{\Lambda}{\nu_{\ll}(\Omega)}\frac{d\mu}{d\tilde{\nu
}_{\ll}}\right)  d\tilde{\nu}_{\ll}\right) \\
&  =\nu_{\perp}(\Omega)u\left(  \infty\right)  +\nu_{\ll}(\Omega)\sup
_{\Lambda>0}\left(  \Lambda+\int_{\Omega}u^{\ast}\left(  \Lambda\frac{d\mu
}{d\tilde{\nu}_{\ll}}\right)  d\tilde{\nu}_{\ll}\right) \\
&  =\nu_{\perp}(\Omega)u\left(  \infty\right)  +\nu_{\ll}(\Omega)\Delta
_{u}(\mu,\tilde{\nu}_{\ll})\text{ .}%
\end{align*}
\end{proof}

\begin{proposition}
\label{Fri}If $u$ has reasonable asymptotic elasticity, then for any $\mu
,\nu\in M^{1}$ such that $\mu\ll\nu$%
\begin{align*}
\Delta_{u}(\mu,\nu)  &  =N_{u}\left(  \nu\parallel\mu\right)  \text{ ,}\\
\delta_{u}\left(  \nu,\mu\right)   &  =e^{H_{u}\left(  \nu\parallel\mu\right)
}-1\text{ .}%
\end{align*}
\end{proposition}

\begin{proof}
First we shall prove the proposition in the case when $\mu,\nu\in M^{1}$ are
equivalent measures. Let $f=\frac{d\nu}{d\mu}$. By Proposition~\ref{ordrel}
above and by Lemma~17 and Theorem~20.4 in~\cite{SloZas04}, we then have
\begin{align*}
N_{u}\left(  \nu\parallel\mu\right)   &  =N_{u}\left(  f\right)
=\sup_{\Lambda>0}\left(  \Lambda+\int_{\Omega}u^{\ast}\left(  \Lambda
/f\right)  fd\mu\right) \\
&  =\sup_{\Lambda>0}\left(  \Lambda+\int_{\Omega}u^{\ast}\left(  \Lambda
\frac{d\mu}{d\nu}\right)  d\nu\right)  =\Delta_{u}(\mu,\nu)\text{ .}%
\end{align*}
Now for any $\mu,\nu\in M^{1}$ such that $\mu\ll\nu$ we take the Lebesgue
decomposition $\nu=\nu_{\perp}+\nu_{\ll}$ into the singular part $\nu_{\perp}
$ and absolutely continuous part $\nu_{\ll}$ with respect to~$\mu$. Then $\mu$
and $\nu_{\ll}$ are equivalent measures. It follows by
Propositions~\ref{relord} and~\ref{Deltadecomp} that
\begin{align*}
N_{u}\left(  \nu\parallel\mu\right)   &  =\nu_{\perp}\left(  \Omega\right)
u\left(  \infty\right)  +\nu_{\ll}\left(  \Omega\right)  N_{u}\left(
\frac{\nu_{\ll}}{\nu_{\ll}\left(  \Omega\right)  }\parallel\mu\right) \\
&  =\nu_{\perp}\left(  \Omega\right)  u\left(  \infty\right)  +\nu_{\ll
}\left(  \Omega\right)  \Delta_{u}\left(  \mu,\frac{\nu_{\ll}}{\nu_{\ll
}\left(  \Omega\right)  }\right)  =\Delta_{u}\left(  \mu,\nu\right)  \text{ .}%
\end{align*}
The equality $\delta_{u}\left(  \nu,\mu\right)  =e^{H_{u}\left(  \nu
\parallel\mu\right)  }-1$ now follows immediately from the definitions of
$\delta_{u}\left(  \nu,\mu\right)  $ and $H_{u}\left(  \nu\parallel\mu\right)
$.
\end{proof}

\section{Concluding remarks}

The notion of $u$-entropy of a probability density, based on the concept of
expected utility maximisation in finance, was first introduced in
\cite{SloZas04} and linked with the Second Law of thermodynamics. In this
paper the definition of $u$-entropy has been extended, on the one hand, to the
case of relative $u$-entropy of one probability measure with respect to
another, and, on the other hand, in the discrete case, to absolute $u$-entropy
of a probability measure. Having established the basic properties of these
notions, we have studied the relationships with other entropy-like quantities
of a similar kind that can be found in the literature. In particular, although
all these approaches yield the Boltzmann-Shannon entropy when the logarithmic
utility is used, it is only the relative $u$-entropy introduced in
Definition~\ref{defreluentr} that is consistent with the R\'{e}nyi entropy for
isoelastic utility functions. The relationships between the various approaches
are summarized in the diagram below. In this context, relative $u$-entropy
emerges as the general unifying quantity among the various approaches related
to expected utility maximisation. \bigskip

\noindent$\!\!\!\!\!
\begin{tabular}
[c]{cccccc}%
$%
\begin{array}
[c]{c}%
\text{general}\\
\text{case}%
\end{array}
\!\!\!$ & $H_{u}\left(  f\right)  $ &
\begin{tabular}
[c]{l}%
$\underrightarrow{\text{Thm.\ref{relord}}}$\\
$\overleftarrow{\,\text{Prop.\ref{ordrel} }}$%
\end{tabular}
& $H_{u}\left(  \nu\parallel\mu\right)  $ & $%
\begin{array}
[c]{c}%
\underrightarrow{\text{ \ Prop.\ref{Fri} \ \ }}\\
\text{\phantom{P}}%
\end{array}
$ & $\delta_{u}\left(  \nu,\mu\right)  $\\
&  &  & \multicolumn{1}{r}{$\left\downarrow \phantom{\displaystyle\int
}\right.  $ \ \ } &  & \\
&  &  & $H_{u}\left(  p\parallel q\right)  $ & $%
\begin{array}
[c]{c}%
\underrightarrow{\text{Prop.\ref{SloZasFri}}(1)}\\
\text{\phantom{P}}%
\end{array}
$ & $D_{u}\left(  p\parallel q\right)  $\\
$%
\begin{array}
[c]{c}%
\text{discrete}\\
\text{case}%
\end{array}
\!\!\!$ &  & \multicolumn{2}{r}{$\text{Prop.\ref{discon}}\left\downarrow
\phantom{\displaystyle\int}\right.  $ \ \ } &
\multicolumn{2}{r}{\cite{FriHuaSan07}$\left\downarrow \phantom{\displaystyle
\int}\right.  $ \ \ }\\
& $H_{-u}^{A}\left(  p\right)  $ & $%
\begin{array}
[c]{c}%
\underleftarrow{\text{Prop.\ref{Ari}}}\\
\text{\phantom{P}}%
\end{array}
$ & $h_{u}\left(  p\right)  $ & $%
\begin{array}
[c]{c}%
\underrightarrow{\text{Prop.\ref{SloZasFri}}(2)}\\
\text{\phantom{P}}%
\end{array}
$ & $H_{u}\left(  p\right)  $%
\end{tabular}
$

\bibliographystyle{amsalpha}
\bibliography{EntrFin}
\end{document}